\RequirePackage{ifpdf}
\ifpdf 
\documentclass[pdftex]{sigma}
\else
\documentclass{sigma}
\fi

\def\cA{\mathcal{A}}            
    \def\cS{\mathcal{S}}

         \def\la{\lambda}                
 
\DeclareMathOperator{\diag}{diag} 

\newcommand{\dsl}[1]{{{\displaystyle{#1}}}}

\def\one#1{#1^{\raise5pt\hbox{$\scriptstyle\!\!\!\!1$}}\,{}}
\def\two#1{#1^{\raise5pt\hbox{$\scriptstyle\!\!\!\!2$}}\,{}}
\def\MMs{Manin\ matrices}
\def\MM{Manin\ matrix}

\def\bea{\begin{eqnarray}}
\def\eea{\end{eqnarray}}
\newcommand{\CC}{{{\mathbb{C}}}}

\newcommand\alg{{{\mathfrak{g}}}}

\def\g{\mathfrak{gl}_n}

\def\gg{\alg}
\def\A{{\mathcal{A}}}
\def\cc{{\mathbb{C}}}

\def\MMs{Manin\ matrices}
\def\MM{Manin\ matrix}

\numberwithin{equation}{section}

\begin{document}

\allowdisplaybreaks

\renewcommand{\thefootnote}{$\star$}

\renewcommand{\PaperNumber}{029}

\FirstPageHeading

\ShortArticleName{Limits of Gaudin Systems: Classical and Quantum Cases}

\ArticleName{Limits of Gaudin Systems:\\ Classical and Quantum Cases\footnote{This paper is a contribution to the Proceedings of the XVIIth International Colloquium on Integrable Systems and Quantum Symmetries (June 19--22, 2008, Prague, Czech Republic). The full collection
is available at
\href{http://www.emis.de/journals/SIGMA/ISQS2008.html}{http://www.emis.de/journals/SIGMA/ISQS2008.html}}}

\Author{Alexander CHERVOV~${}^\dag$, Gregorio FALQUI~${}^\ddag$ and Leonid RYBNIKOV~${}^\dag$}

\AuthorNameForHeading{A. Chervov, G. Falqui and  L. Rybnikov}

\Address{$^\dag$~Institute for Theoretical and Experimental Physics, \\
\hphantom{$^\dag$}~25 Bolshaya   Cheremushkinskaya Str., 117218 Moscow, Russia}

\EmailD{\href{mailto:chervov@itep.ru}{chervov@itep.ru}, \href{mailto:leo.rybnikov@gmail.com}{leo.rybnikov@gmail.com}}

\Address{$^\ddag$~Dipartimento di Matematica e Applicazioni,
Universit\`a di Milano -- Bicocca,\\
\hphantom{$^\ddag$}~via R. Cozzi, 53, 20125 Milano, Italy}
\EmailD{\href{mailto:gregorio.falqui@unimib.it}{gregorio.falqui@unimib.it}}
\URLaddressD{\url{http://www.matapp.unimib.it/~falqui/}}

\ArticleDates{Received November 01, 2008, in f\/inal form February 25,
2009; Published online March 09, 2009}

\Abstract{We consider the XXX homogeneous Gaudin system with $N$
sites, both in classical and the quantum case. In particular we show
that a suitable limiting procedure for letting the poles of its Lax
matrix collide can be used to def\/ine new families of Liouville
integrals (in the classical case) and new ``Gaudin'' algebras (in
the quantum case). We will especially treat the case of total
collisions, that gives rise to (a generalization of) the so called
Bending f\/lows of Kapovich and Millson. Some aspects of multi-Poisson
geometry will be addressed (in the classical case). We will make use
of properties of ``Manin matrices'' to provide explicit generators
of the Gaudin Algebras in the quantum case.}

\Keywords{Gaudin models; Hamiltonian structures; Gaudin algebras}

\Classification{82B23; 81R12; 17B80; 81R50}

\renewcommand{\thefootnote}{\arabic{footnote}}
\setcounter{footnote}{0}

\section{Introduction}
The Gaudin model \cite{G83} was introduced by M.~Gaudin as a spin model related
to the Lie algeb\-ra~$sl_2$, and later generalised to the case of an
arbitrary semisimple Lie algebra $\alg$.

The Hamiltonian is
\begin{equation*}
H_G= \sum_{a=1}^{\dim  \alg}\sum_{i \neq j} x_a^{(i)} {x^a}^{(j)},
\end{equation*}
where $\{x_a\}$, $a=1,\dots,\dim\gg$, is an orthonormal basis of
$\gg$ with respect to the Killing form (and~$x^a$ its dual). These
objects are regarded as elements of the polynomial algebra
$\cS(\alg^*)^{\otimes\, N}$ in the classical case, and  as elements
of the universal enveloping algebra $U(\alg)^{\otimes N}$ in the
quantum case, as
\begin{equation*}
\dsl{ x_a^{(i)}=1\otimes\cdots\otimes \!\!\!\underbrace{x_a}_{i\text{-th factor}}\!\!\!\otimes\, 1\cdots\otimes 1.}
\end{equation*}

Gaudin himself found that the quadratic Hamiltonians
\begin{equation*}
H_i=\sum\limits_{k\neq i}\sum\limits_{a=1}^{\dim\gg}
\frac{x_a^{(i)}{x^a}^{(k)}}{z_i-z_k}
\end{equation*}
provide a set of
quantities that commute with $H_G$, for any choice of pairwise
distinct points  $z_1,\dots,z_N$ in the complex plane. The
corresponding classical quantities are  ``constants of the motion''
for the Hamiltonian f\/low generated by $H_G$.
For instance, in the $sl_2$ case, we have
\begin{gather*}
 H_G=\sum_{i,j=1,\, i\neq j}^N h^{(i)}h^{(j)}+e^{(i)}f^{(j)}+f^{(i)}e^{(j)},\nonumber   \\
 H_i=\sum\limits_{k=1, k\neq
   i}^{N}\frac{h^{(i)}h^{(k)}+e^{(i)}f^{(k)}+f^{(i)}e^{(k)}}{z_i-z_k},
   \qquad i=1,\ldots,N. 
\end{gather*}

Later it was shown (see, e.g., \cite{Ju89}) that -- in the classical case -- for a
general (semisimple) Lie algebra $\alg$ the spectral invariants of
the Lax matrix
\begin{equation}\label{LG}
 L_G(z)=\sum_{i, a} \frac{x_a^{(i)}\otimes {x^a}^{(i)}}{z-z_i}
\end{equation}
encode a (basically complete) set of invariant quantities. We recall
that under the name spectral invariants of $L_G$ we mean the
quantities
\begin{equation*}
\mathop{\rm Res}\limits_{z=z_i} (z-z_i)^{k_\alpha}
\Phi_\alpha(L_G),\qquad i=1,\ldots, N,
\end{equation*} where
$\Phi_\alpha$ is a complete set of Ad-invariant quantities for
$\alg$, and the exponents $k_\alpha$ run in an appropriate set. For
instance, for $\alg=gl(r)$, which is the case we will basically
consider, we can choose $\Phi^\alpha(L_G(z))$ to be
$\text{Tr}\,(L_G(z)^\alpha)$ (or, equivalently, the coef\/f\/icients
$p_\alpha$ of the characteristic polynomial of $L_G(z)$), with
$\alpha=1,2,\ldots, r$; in this case, $k_\alpha$ will run from $0$
to $\alpha-1$.

The meaning of ``basically'' complete refers to the fact that, e.g.,
the residues of the trace of~$L^2_G(z)$ at the dif\/ferent points~$z_i$,
$i=1,\ldots, N$ are not independent, since $L_G(z)$ is
regular at inf\/inity. Actually, the spectral invariants of $L_G(z)$
generate the Poisson-commutative subalgebra of maximal possible
transcendence degree in the algebra of diagonal invariants
$S(\alg^*\oplus\cdots\oplus\alg^*)^\alg$. To obtain a maximal
Poisson-commutative subalgebra in
$S(\alg^*\oplus\cdots\oplus\alg^*)$, one has to add
$\frac12(\dim\g+{\rm rank}\,\g)$ independent Hamiltonians (the
Mischenko--Fomenko generators for the diagonal $\alg$).
This is due to the fact that, being the models we are considering
homogeneous, they are invariant with respect to the diagonal adjoint
action of the group $G$. We will later refer to this fact as
the {\em global} $G$ invariance of the models.
The enlarged set herewith def\/ined gives to the (classical) Gaudin
system the structure of a Liouville integrable system. For instance,
in the case $\alg=sl_2$, one has to complement the residues of the
trace of $L_G^2$ at its singular points~$z_i$ with a further
quantity, e.g.\ a component of the total ``spin'', say
$S_z=\sum\limits_{i=1}^N h^{(i)}$. It should however be noticed that the
number of integrals to be added to the spectral invariants equals,
independently of the number $N$ of ``Gaudin generalized magnets'', the
rank of~$\alg$. So, with a slight abuse of language, we will
sometimes speak of a complete set of integrals referring only to the
spectral invariants.

The quantum case was the subject of intense investigations. It is
outside of the size of this paper to historically frame such a
research line. For our purposes, it is of paramount relevance the
fact that, in \cite{FFR94}, Feigin, Frenkel and Reshetikhin proved
the existence of a large commutative subalgebra
$\cA(z_1,\dots,z_N)\subset U(\gg)^{\otimes N}$ containing the
quadratic elements $H_i$  (see, also, \cite{ER96, FF,FFTL,Fr95,Ryb1,Ryb}).
 For $\gg=sl_2$, the
algebra $\cA(z_1,\dots,z_N)$ is generated by the $H_i$, the
additional global element $S_z$ and the central elements of
$U(\gg)^{\otimes N}$.

In the other cases, the algebra $\cA(z_1,\dots,z_N)$ has also some
new generators known as higher Gaudin Hamiltonians. Their explicit
construction for $\gg=gl(r)$ was obtained by D.~Talalaev~\cite{Ta04}
and further discussed in papers by A.~Chervov and D.~Talalaev~\cite{CT04,CT06}.

In the present paper we will
discuss
the problem of considering -- both from the classical and from the
quantum point of view -- what happens when the arbitrary points
$z_1,\ldots,z_N$ appearing in the Lax matrix, and in the quadratic
Hamiltonians $H_i$ glue together; in particular, we will pay special
attention to  the ``extreme'' case, when in some sense all points
collide; indeed, this extreme case gives rise to (the generalisation
of) a remarkable integrable system, called the ``bending f\/lows''
systems introduced and studied~-- for the case $\alg=so(3)$~-- by
Kapovich and Millson~\cite{KM96} in their work concerning moduli
spaces of polygons in ${\mathbb R}^3$.

\section{The classical case}\label{sect-class}

In the classical case, to shorten notations, it is customary to write the
Lax matrix $L_G(z)$ \eqref{LG} as $L_G(z)= \sum\limits_{i=1}^N
\frac{X_i}{z-z_i}$.  It satisf\/ies the Poisson algebra
\begin{equation*}
\{L_G(z)\otimes 1, 1 \otimes L_G(u)\} =\left[\frac{\Pi}{z-u}, L_G(z)\otimes 1 +
1 \otimes L_G(u)\right] ,
\end{equation*}
where $\Pi$ is the  permutation matrix $\Pi(X\otimes Y)=Y\otimes
X$.

As it is well known, this linear $r$-matrix structure (together with
suitable reductions) is associated with a huge number of classical
integrable systems, such as Neumann type systems and the $n$-dimensional
Manakov tops, f\/inite gap reductions of
the KdV equations, as well f\/inite gap reductions of its
generalisations to arbitrary simple Lie algebras known as
Gel'fand--Dickey hierarchies,
Hitchin's system on singular rational curves and so on and so forth
(see, e.g., \cite{RSTS,BBT03} and the references quoted therein).

Also, it is well known that the $r$-matrix Poisson brackets
presented above can be seen as a~kind of ``shorthand notation'' for
the following situation:
\begin{itemize}\itemsep=0pt
\item
The phase space\footnote{In this paper we will always deal with
metric (or reductive) Lie
  algebras (and, in particular, with $gl(r)$), so that we will
  tacitly henceforth identify $\alg$ with $\alg^*$.}
 of the $N$-site classical Gaudin model is ${\alg^*}^{\otimes\,N}$.

\item The ``physical'' Hamiltonian is a mean f\/ield spin-spin interaction,
\begin{equation*}
H_G=\frac12\sum_{i\neq j=1}^N \text{Tr}\,(X_i\cdot X_j),\qquad  X_j\in
\alg^* (\simeq \alg).
\end{equation*}
\item
The Poisson brackets def\/ined by the $r$-matrix formula are just
product of Lie--Poisson brackets on ${\alg^*}^{\otimes\,N}$.
\item
The def\/inition of the (classical) Lax matrix as $ L_G(z)=\sum\limits_{i=1}^N
\frac{X_i}{z-z_i}$ def\/ines an embedding of our phase space
${\alg}^{\otimes\,N}$ into a Loop space $L\alg\simeq\alg((z))$.
\end{itemize}
We now turn to the discussion of  limits of the above situation when
some of the points $z_1,\ldots,z_N$ collide (or glue together) in a
suitable sense. We shall discuss both algebraic aspects (that is,
the corresponding Lax matrices) and address also some Poisson
aspects.

\subsection{Lax matrices and their limits}\label{subs-laxlim}

\begin{figure}[t]
\centerline{\includegraphics[width=5cm]{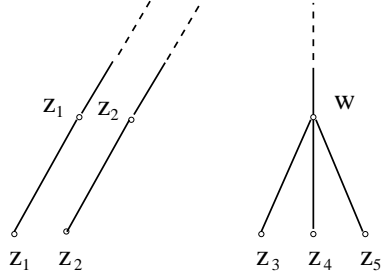}}

  \caption{A gluing pattern $\pi$.}\label{figurgr}
\end{figure}

To consider limits of the Gaudin algebras when some of the points
$z_1,\dots,z_N$ glue together, we can proceed as follows. We keep
some (say the f\/irst $k$) points $z_1,\dots,z_k$ ``f\/ixed'', and let the
remaining $N-k$ points glue to a new point $w$, via{\samepage
\begin{equation}\label{gluing} z_{k+i}=w+s\,u_i,\quad i=1,\dots,N-k,\quad z_i\neq z_j, \quad u_i\neq u_j,  \quad
s\to 0 .
\end{equation}
We can represent this procedure as in Fig.~\ref{figurgr}.}

To see what happens in this limit we f\/irst observe,
that obviously
\begin{equation}\label{lax-gluing1}
L_{G}(z) \to L_2(z)= \sum_{i=1}^{k}\frac{X_i}{z-z_i}
+\frac{\sum\limits_{i=k+1}^N X_i}{z-w}, \qquad s\to 0. \end{equation}

Considering only this limit would however be too {\em na\"ive}, and indeed, we see that (even in the case of $\alg=gl(2)$) the number of
Hamiltonians obtained from the spectral invariants of \eqref{lax-gluing1} is not suf\/f\/icient to yield complete integrability.

One can notice that, in the gluing procedure, some residues of the Lax matrix should also be taken into account.
That this is the case can be seen by a proper rescaling argument as follows.

Let us introduce a new variable $\tilde z = (z-w)/s$,
so that $z=w+s \tilde z$, and rewrite the Lax matrix in terms of the new
variable as
\begin{gather*}
 L_{G}(z) = \sum_{i=1}^{k}\frac{X_i}{w +s \tilde z-z_i}
+\sum_{i=k+1}^N\frac{ X_i}{w+s \tilde z -w -su_i}.
\end{gather*}
Multiplying this by  $s$ we can see that another Lax matrix appears
for $s\to0$, that is,
\begin{gather*}
 sL_{G}(\tilde z)  = s
\left(\sum_{i=1}^{k}\frac{X_i}{w +s \tilde z-z_i} + \sum_{i=k+1}^N\frac{
X_i}{s(\tilde z -u_i) }\right) \to \sum_{i=k+1}^N \frac{X_i}{\tilde z -u_i
} =L_1(\tilde z),
\qquad s\to 0.
\end{gather*}
The justif\/ication of this rescaling and limit
is simple.
Indeed,  for $s\neq 0$, the ring of Hamiltonians obtained from
$sL_{G}(\tilde z)$ coincides with that obtained from
$L_{G}( z)$. Thus the Hamiltonians
obtained from the limiting Lax matrix $L_1(\tilde z)$
are limits of the corresponding Hamiltonians of the
Gaudin model.

Finally we conclude that starting from the Lax matrix with
generic (distinct) points $z_1,{\dots},z_N,$  we  can associated, to the
gluing of \eqref{gluing} the following pair of ``Lax matrices'':
\begin{equation}
\label{lax-gluing}
L_{1}(z)=\sum_{i=k+1}^N \frac{X_i}{z-u_i}\qquad\text{and} \qquad
\quad L_2(z)=\sum_{i=1}^{k}\frac{X_i}{z-z_i} +\frac{\sum\limits_{i=k+1}^N
X_i}{z-w}.
\end{equation}
Notice that, in the general case, we can choose the gluing procedure
to be explicitly given by, e.g.,
\begin{equation*}
z_{k+i}=w+s(z_{k+i}-w),\qquad s\in(0,1)
\end{equation*}
and, using invariance w.r.t.\ transformation of the spectral
parameter $z\to z-w$, trade the matrix $L_1$ of \eqref{lax-gluing}
for $ \tilde{L}_{1}(z)=\sum\limits_{i=k+1}^N \frac{X_i}{z-z_i}$.

In particular, in the example of the Fig.~\ref{figurgr}, we
would associate, to the ``generic'' Lax matrix $ L_G(z)=\sum\limits_{i=1}^5
\frac{X_i}{z-z_i} $ the two matrices
\begin{gather*}
L_1(z)=\frac{X_3}{z-z_3}+\frac{X_4}{z-z_4}+\frac{X_5}{z-z_5},\qquad \text{and} \qquad
L_2(z)=\frac{X_1}{z-z_1}+\frac{X_2}{z-z_2}+\frac{X_3\!+X_4\!+X_5}{z-w}.
\end{gather*}

\begin{proposition}\label{prop1}
For every choice of $w\in \CC$ the family of
spectral invariants $\mathcal{H}^{(1)}$, $\mathcal{H}^{(2)}$ associated
with the Lax matrices $L_1$
and $L_2$
satisfy the following properties:
\begin{enumerate}\itemsep=0pt
\item[$1.$]
The elements of $\mathcal{H}^{(1)}$, $\mathcal{H}^{(2)}$ commute
w.r.t.\ the standard (diagonal) Poisson brackets on~$\alg^N$.

\item[$2.$] The  dimension  of the Poisson commutative subalgebra
${\mathcal H}_{1,2,w}$ generated by the spectral invariants
  $\mathcal{H}^{(1)}$ and $\mathcal{H}^{(2)}$ which are those, respectively, associated  with the Lax matri\-ces~$L_1(z)$
  and $L_2(z)$, coincides with that of the
spectral invariants associated with the generic Lax matrix.

\item[$3.$] The physical Hamiltonian $H_G$ lies in~${\mathcal H}_{1,2,w}$.
\end{enumerate}
\end{proposition}

\begin{proof}
The commutativity property (the f\/irst assertion) holds thanks to the
following facts. First, it is clear that if $h_1$, $h_2$ are taken
either both from the subfamily~$\mathcal{H}^{(1)}$ or the subfamily~$\mathcal{H}^{(2)}$ the assertion is trivially true. Let us suppose
thus that $h_1\in \mathcal{H}^{(1)}$ and $h_2\in \mathcal{H}^{(2)}$.
The fact that $\{h_1,h_2\}=0$ follows from the fact that $h_2$
depends on $X_1,X_2,\ldots,X_k$ only through the sum
$\mathbf{X}_k=\sum\limits_{i=1}^k{X_i}$, and $h_1$ is invariant under the
diagonal action of~$\alg$.

The second assertion holds thanks to the results of Jur\v{c}o brief\/ly
mentioned before about the number of functionally independent
spectral invariants of a Lax matrix of Gaudin type, as well as to
the functional dependence of $L_2(w)$ on the ``variables''
$X_{1},\ldots,X_{k}$.

The third one can be checked by induction. \end{proof}

{\bf Remarks}. 1) The meaning of this proposition is that, once the
global $G$-invariance has been taken into account, the procedure of
glueing discussed above provides new families of Liouville integrals
for the ``physical'' Gaudin Hamiltonian $H_G$.

2) It is fair to say that a procedure somewhat similar to this was presented
by the late V.~Kuznetsov in~\cite{Ku92}. However, possibly owing to the fact
that the author restricted himself to the case of rank~$1$ (namely, mostly
to the case $\alg=so(2,1)$), he did not consider the need to add the further
singular point $w$. This is crucial to obtain, in the general case, the
correct number of Liouville integrals for the ``new'' system.

3) Our glueing procedure dif\/fers substantially from the one
considered in~\cite{MPR05}
(see also~\cite{Ch03}), where the limiting procedure produces higher
order poles in the Lax matrix, and is associated with non-semisimple
Lie--Poisson algebras, that are In\"on\"u--Wigner contractions of the
original one (see, for our case, Section~\ref{subs-poa}).

\subsection{On the Poisson geometry of the limiting procedure}\label{subs-poa}

The interplay between classical Lax matrices, Loop algebras and
Poisson manifolds is nowadays well known, and was settled mainly in
the  works of the Leningrad's school (see, e.g.\ the review by Reyman
and Semenov-Tyan-Shanski in~\cite{RSTS}). It is encoded in the notion of
 $R$-operator as follows.

On the space ${\mathfrak{g}}^*((z))$ of Laurent polynomials with
values in (the dual of) a Lie algebra $\alg$, there is a family of
mutually compatible Poisson brackets, $\{\cdot,\cdot\}_k$ associated
with a family of classical $R$-operators
\begin{equation*}
R_k(X(z))=\big(z^k X(z)\big)_+-\big(z^k X(z)\big)_-
\end{equation*}
via the formula
\[
\{F,G\}_k(X)=\langle X,[R_k(\nabla(F)),\nabla(G)]\rangle- \langle
X,[R_k(\nabla(G)),\nabla(F)]\rangle.
\]
Spectral invariants of a Lax matrix (say, polynomially dependent on
the spectral parameter $z$) form an Abelian Poisson subalgebra
w.r.t.\ any of the brackets $\{\cdot,\cdot\}_k$.

Let us consider the space of $N$-th order polynomials
$\mathfrak{g}^*_N= \alg^* [[z]]\text{ mod } z^{N+1}$. It can be
shown that (see, e.g., \cite{RSTS,PV})
\begin{itemize}\itemsep=0pt
\item The brackets associated with
$R_0,\ldots,R_{N}$ on ${\mathfrak{g}}^*((z))$
 restrict to the af\/f\/ine subspace
\begin{equation*}
\alg^*_{N,A}:=\left\{X\in \alg^*((z))\,  |\, X(z)=z^{N+1}
A+\sum_{i=0}^{N}z^i X_i\right\},
\end{equation*}
where $A$ is a {\em fixed} element of $\mathfrak{g}$, and thus on
$\alg^*_N=\alg^*_{N,\mathbf{0}}$, that is when $A=0$ which is the
case we plan to consider in this paper.
\item
These brackets are mutually compatible and give rise to
multi-Hamiltonian structures on the f\/inite dimensional manifolds
$\alg^*_{N,A}$. A straightforward observation~\cite{FM03-2} is the
following: we can associate with any polynomial of degree $N$
${\cal Q}(z)=\sum\limits_{i=0}^N \kappa_i z^{N-i}$ a Poisson bracket (of
the RSTS family) via:
\begin{equation*}
\{\cdot,\cdot\}_{\cal Q}=\sum_{i=0}^N \kappa_i\{\cdot , \cdot\}_i.
\end{equation*}
\end{itemize}
In this setting, we can recover  the standard (product) Lie--Poisson
structure on $\alg^{\otimes\,N}$ as follows. Considering the
``Lax map'' $L_G=\sum\limits_{i=1}^N \frac{X_i}{z-z_i}$
the diagonal structure can be obtained, in the RSTS
framework, by the sum
\begin{equation*}
\{\cdot,\cdot\}_{\cal S}=\{\cdot,\cdot\}_{N}+\sum_{i=0}^{N-1}
(-1)^i\sigma_i\{\cdot , \cdot\}_i,
\end{equation*}
the $\sigma_i$ being the elementary symmetric polynomials in the
quantities $z_1,\ldots,z_N$. That is, the standard $r$-matrix
Poisson bracket can be regarded as the one associated with the
polynomial
\begin{equation*} S(z)=\prod_{i=1}^N
(z-z_i)=z^N+\sum_{i=0}^{N-1} (-1)^i\sigma_i z^{N-1}.
\end{equation*}
As a side remark, one can also notice that this gives the
possibility of constructing a bi-Hamiltonian structure for the
classical Gaudin case. Indeed the structure def\/ined by the
polynomial $\left(\frac{S(z)}{z}\right)_+$ provides, together
with the standard one, a bi-Hamiltonian structure for the Gaudin
model, such that the spectral invariants of $L_G(z)$ f\/ill in
recursion relations of (genera\-li\-zed) Lenard--Magri type. More in
detail, suitable combinations thereof give rise to {\em anchored}
Lenard--Magri sequences, in the terminology of
Gel'fand--Zakharevich~\cite{GZ00}.

What is more important for the present paper is that this picture
suggests and gives the opportunity of studying, from this point of
view, the limits of (suitable combinations of) the RSTS Poisson
structures when some (and, iteratively, possibly all) poles of the
Lax matrix $z_1,\ldots,z_N$ glue together.

This means the following: if we pull back, via the map
\begin{equation*}
L_G\to L_{G,{\rm pol}},\qquad L_{G,{\rm pol}}=\prod_{i=1}^N
(z-z_i)\left(\sum_{i=1}^N \frac{X_i}{z-z_i}\right)
\end{equation*}
the RSTS Poisson brackets $\{\cdot,\cdot \}_{k}$, $k=0, \ldots ,N$, we
obtain families of {\em linear} Poisson brackets on
$\mathfrak{g}^*_N= \alg^* [[z]]\ \text{mod} \ z^{N+1}$  that depend
rationally on $(z_1,\ldots,z_N)$.

In this framework, a natural question arises, that is, whether one
can associate suitable Poisson structures to the procedure of gluing
(some of the) points $z_i$ discussed -- in the Lax setting -- in
Section~\ref{subs-laxlim}. In the example of the limit pattern
depicted in Fig.~\ref{figurgr}, one of the possible Poisson
structures\footnote{There is two-parameter family of such
structures.} that can be obtained is represented by the operator (for
the sake of simplicity, we set $w=z_3$, $\alg=gl(r)$)
\begin{equation}\label{neq0}
\mathcal{P}= \left[ \begin {array}{ccccc}
0&z_{{23}}[X_{{1}},\cdot]&0&0&0
\\
\noalign{\medskip}z_{{23}}[X_{{1}},
\cdot]&
[P_{22},\cdot]&-z_{{12}}[X_{{3}},\cdot]&-z_{{12}}[X_{{4}},\cdot]
&-z_{{12}}[X_{{5}},\cdot]
\\
\noalign{\medskip}0&-z_{{12}}[X_{{3}},\cdot]&0&-\dsl{{\frac
{z_{{13}}z_{{34} }}{z_{{45}}}}}[X_{{3}},\cdot]&\dsl{{\frac
{z_{{13}}z_{{35}}}{z_{{45}}}}}[X_{{3}},\cdot]
\\
\noalign{\medskip}0&-z_{{12}}[X_{{4}},\cdot]&-\dsl{{\frac
{z_{{13}}z_{{34}
}}{z_{{45}}}}}[X_{{3}},\cdot]&[P_{44},\cdot]&[P_{45},\cdot]
\\
\noalign{\medskip}0&-z_{{12}}[X_{{5}},\cdot]&\dsl{{\frac
{z_{{13}}z_{{35}}}{z_{{45}}}}}[X_{{3}},\cdot]&[P_{45},\cdot]&[P_{55},\cdot]\end
{array} \right],
\end{equation}
where $z_{ij}=z_i-z_j$, and
\begin{gather*}
P_{22}=z_{{23}}\left(X_{{2}}-X_{{1}}\right)+z_{{12}}\left(X_3+X_4+X_5\right),\nonumber\\
P_{44} = {\frac{z_{{13}}z_{{34}}}{z_{{45}}}}\left(X_3-\frac{
z_{{35}}-z_{{45}}}{{z_{{45}}}} X_4-{\frac
{z_{{13}}{z_{{34}}}}{{z_{{45} }}}}X_5\right),\nonumber
\\
P_{45}=\frac{z_{{13}}}{{{{z^2_{{45}}}}}} \big(
z_{{35}}^{2}X_{{4}}+ z_{{34}}^{2}X_{{5}}\big),\\
P_{55}=-{\frac{z_{{13}}z_{{35}}}{z_{{45}} }}\left(X_3+{\frac
{{z_{{35}}}}{{z_{{45}}}}}X_4+{\frac{z_{{34}}-z_{{45}}}{{z_{{45}}}
}}X_{ {5}}\right).\nonumber
\end{gather*}
Using the fact that both tangent vectors and one-forms on the phase
space of this f\/ive-site model can be identif\/ied with f\/ive-tuples of
matrices the meaning of the representation~\eqref{neq0} (see also the
subsequent equations (\ref{xx1}), (\ref{xx2}) can be illustrated as follows.

The Hamiltonian vector f\/ield associated via $\mathcal{P}$
to the one-form $(\alpha_1,\ldots, \alpha_5)$ is given by $
\dot{X}_i=\sum\limits_{j=1}^5 [P_{ij},\alpha_j]$. E.g.,
$\dot{X}_1=z_{23}[X_1,\alpha_2]$, and so on and so forth.

It is not dif\/f\/icult to check that the spectral invariants of the matrices
\begin{gather*}
L_1=\frac{X_3}{z-z_3}+\frac{X_4}{z-z_4}+\frac{X_5}{z-z_4},\qquad \text{and}\qquad
L_2=\frac{X_1}{z-z_1}+\frac{X_2}{z-z_2}+\frac{X_1+X_2+X_3}{z-z_3}.
\end{gather*}
do commute with respect to the Poisson brackets $\cal{P}$ def\/ined by \eqref{neq0} (besides, as proven in Proposition~\ref{prop1}, commuting w.r.t.\ the standard Poisson brackets).

\medskip

\noindent
{\bf Conjecture.} For every limit pattern $\pi$ described in
Section~\ref{subs-laxlim} there are suitable li\-near combinations
$\{\cdot,\cdot\}_\pi$ of the Poisson brackets above that remain
regular. Endowing the phase space~${\alg^{*}}^{\otimes\,N}$  with these limiting brackets, one obtains bi-Hamitonian
(or, possibly, multi-Hamiltonian) manifolds. The integrable systems def\/ined,
(according to the Gel'fand--Zakharevich scheme)  by these
multi-Hamiltonian structures admit Lax representations, whose Lax matrices are exactly those constructed in Section~\ref{subs-laxlim}.
In other words,the spectral invariants of the Lax matrices associated with the limiting pattern $\pi$
patterns provide, w.r.t.\ the Poisson pencil formed out of the
diagonal bracket and the bracket $\{\cdot,\cdot\}_\pi$ Lenard--Magri
sequences (possibly in a generalised sense). Also, the brackets
$\{\cdot,\cdot\}_\pi$ are semisimple Lie--Poisson brackets.

\medskip

This conjecture has not been fully proven yet, but we have checked it in a signif\/icant number of examples.
We notice that, in the ``extreme gluing case'', that is, when all
$z_i$ glue according to the pattern $z_2\to z_1$, $z_3\to
z_1$ $(=z_2)$, $z_4\to z_1$ $(=z_2=z_3)$, $\ldots$ we can obtain the following
``limit'' bracket which is independent of any parameter:
\begin{gather*}
 \{ F, G \}_{\rm limit}= \sum_{i,j,k} r_{ijk} \text{Tr} \left(\nabla
F_i\left[ X_k, \nabla G_j\right] \right)\qquad \text{with } \\   r_{ijk}=(k-1)
\delta_{ij} \delta_{jk}-\theta_{(i-k)} \delta_{ij} +\theta_{(j-i)}
\delta_{ik}+\theta_{(i-j)}\delta_{jk} ,
\end{gather*}
where $\nabla F_i$, (resp.~$\nabla G_i$) represents the
dif\/ferential of $F$ (resp.~$G$) w.r.t.\ the $i$-th entry $X_i$.

As it has been shown in \cite{FM04} these Poisson brackets are still
compatible with the standard $r$-matrix ones.

To elucidate the formula above, we notice that, in the $4$-site case we have the following representation~-- by
Poisson operators~-- of the brackets as follows:
\begin{gather}\label{xx1}
P_{\text{standard } r}=\left(\begin{array}{cccc} [X_1,\cdot]&&&\\
& [X_2,\cdot]&&\\
&&[X_3,\cdot]&\\ &&&[X_4,\cdot]
\end{array}\right),
\\
\label{xx2}
P_{\text{limit}}=\left(
\begin{array}{cccc}0
& [X_1,\cdot\>]& [X_1,\cdot\>]&[X_1\cdot\>] \\
 {} [X_1,\cdot\>]& [X_2-X_1,\cdot\>]& [X_2,\cdot\>]& [X_2,\cdot\>]\\
  {}[X_1,\cdot\>]& [X_2,\cdot\>]& [2X_3-X_2-X_1,\cdot\>]& [X_3,\cdot\>]\\
 {}[X_1,\cdot\>]& {}[X_2,\cdot\>]& [X_3,\cdot\>]&\big[3X_4-\sum\limits_{i=1}^3X_i, \cdot\>\big]
\end{array}\right).
\end{gather}
According to the Gel'fand--Zakharevich (or Lenard--Magri) scheme,
$P_{\text{standard } r}$ and $P_{\text{limit}}$ def\/ine the $gl(r)$
(as well as $\alg$) -- valued generalisation of the so called {\em
Bending flows}, introduced in the case corresponding to
$\alg=\mathfrak{gl}(2)$ by Kapovich and Millson~\cite{KM96}, and further discussed in \cite{FlM01}.

\subsection{The Bending system}

Bending f\/lows are def\/ined on the moduli space $\mathcal{M}_{\mathbf
r}$ of $(N)$-gons with f\/ixed sides lengths
$\mathbf{r}=(r_1,\ldots,r_{N})$, $r_i>0$. In~\cite{KM96} it was shown, among
other properties, that:
\begin{itemize}\itemsep=0pt
\item $\mathcal{M}_{\mathbf
r} $ is a smooth $(2N-6)$-dimensional open manifold, whose compactif\/ication
is achieved when some of the side lengths vanish.
\item
$\mathcal{M}_{\mathbf r}$ is endowed with a natural symplectic structure,
since it is a symplectic quotient of products of spheres:
\[
\big\{(e_1,\,e_2,\dots,e_{N})\in\, S^2_{r_1}\times
S^2_{r_2}\times\cdots \times S^2_{r_{N}}\big\}//SO(3),
\]
where the symbol $//$ denotes the symplectic quotient on the null
manifold of the moment map def\/ined by the diagonal $SO(3)$ action on
the spheres~$S_{r_i}$.
\end{itemize}

\begin{figure}[t]
\begin{minipage}[b]{7cm}
\centerline{\includegraphics[width=6cm]{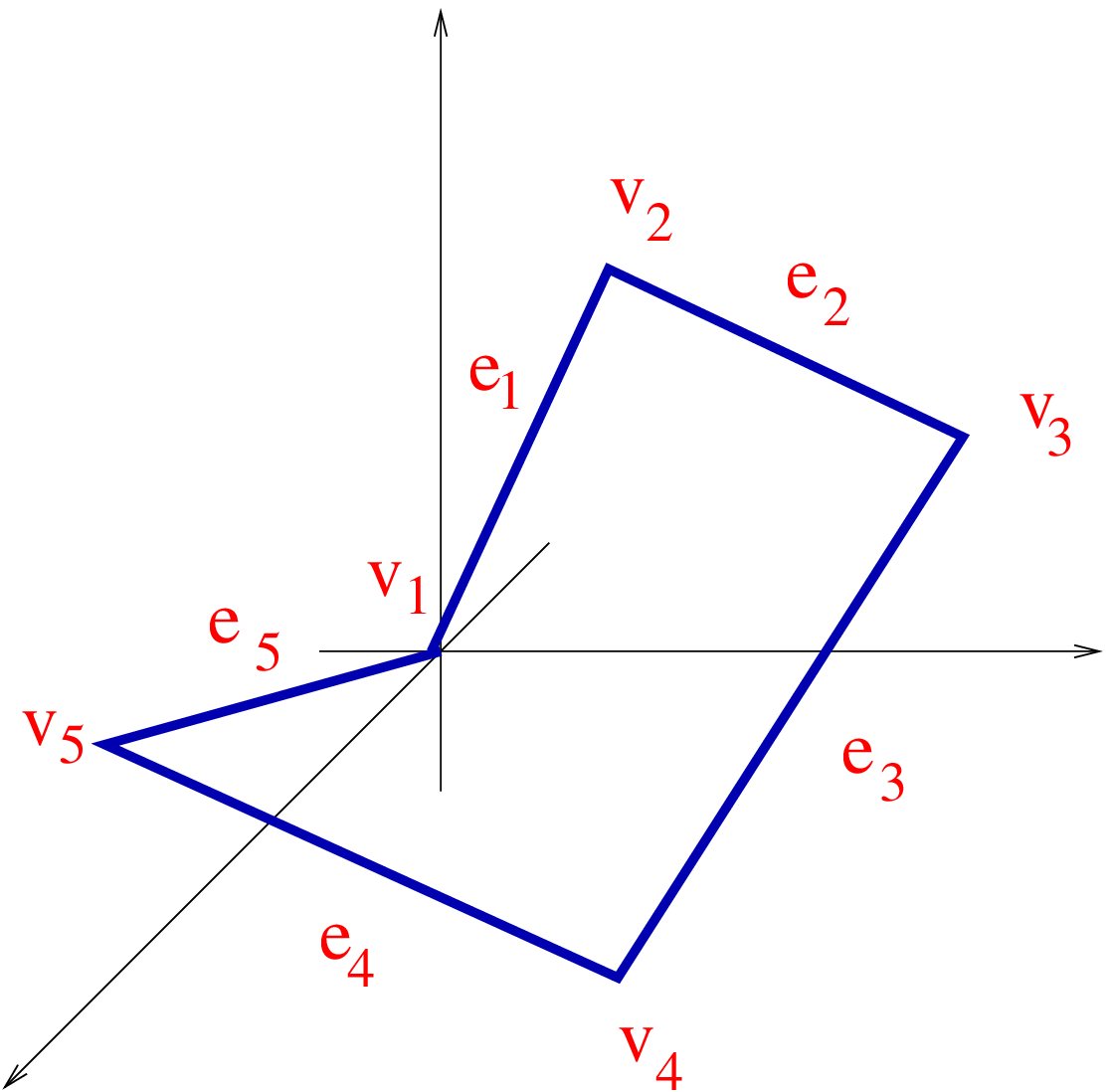}}
\caption{Polygon in $\mathbb{R}^3$.} \label{bspac}
\end{minipage}\,
\begin{minipage}[b]{8.8cm}
\centerline{\includegraphics[width=6cm]{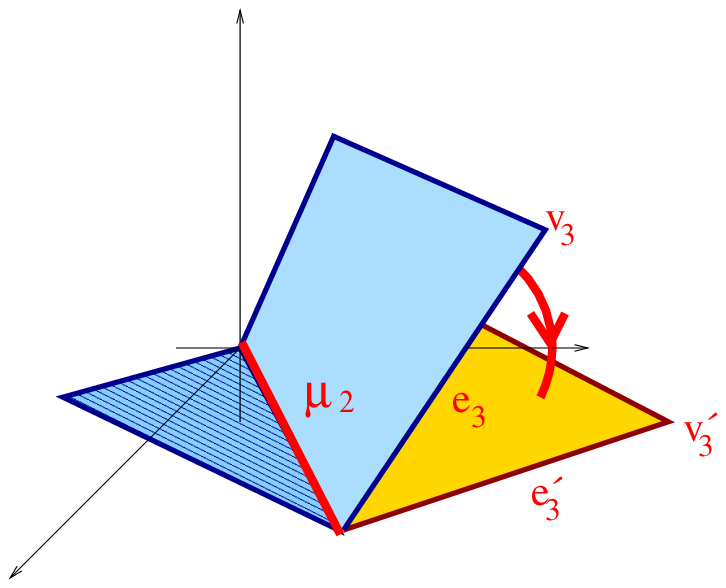}}
\caption{The Bending f\/lows for polygons ($(so(3))$-case).}\label{bflow}
\end{minipage}
\end{figure}

On this moduli space of polygons a completely
integrable Hamiltonian system can be def\/ined, indeed
called the {\em Bending system} or system of bending f\/lows.

The action variables for such system are the lengths of diagonals
stemming from one vertex, and the angle variables are (indeed) the
dihedral angles (see Figs.~\ref{bspac} and~\ref{bflow}). Geometrically, the f\/lows {\em bend} one part of the
polygon around the corresponding diagonal, keeping the rest f\/ixed,
whence the denomination.

%

\begin{figure}[t]

\centerline{\includegraphics[width=10.5cm]{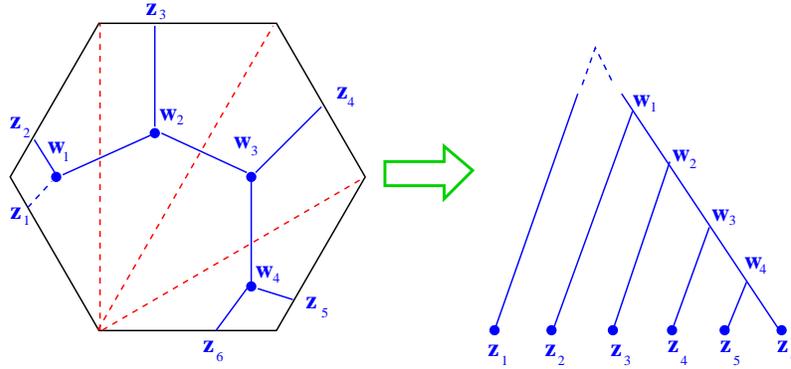}}

 \caption{Triangulations and trees.}\label{fig2.4}
 \end{figure}

In the papers \cite{FM01,FM04} the following picture was
established: $\alg$-Bending f\/lows can be def\/ined on the same phase
space of the $\alg$-Gaudin system, the link between the two (in the
Kapovich Millson $N$-gon case) being the fact that $S^2_{r}$ is a
symplectic leaf of ${so}(3)^N$, and the condition for the
polygon to close is just the choice of the special value $\mu=0$ of
the image of the momentum map. Also, the following results were
established:
\begin{itemize}\itemsep=0pt
\item  $\alg$-Bending  f\/lows admit a set of $N-1$ Lax matrices of
the form
\begin{equation}\label{lax-fm}
 L_k(z)=zX_k+\sum_{i=k+1}^N X_i, \qquad k=1,\ldots,N-1.
 \end{equation}
\item
Hamiltonians come in ``clusters'', each cluster being associated with
the corresponding Lax matrix. For $\alg=gl(r)$
\begin{equation}\label{ham-fm}
H^a_{k,m}=\mathop{\rm res}\limits_{z=0}\frac1{z^{a+1}} \, \text{Tr}\,
L_k^m(z).\end{equation}
\item
Separation of variables can be performed in this scheme. Separation
canonical conjugated variables come in ``clusters'', associated, via
the bi-Hamiltonian version~\cite{FaPe03} of the Sklyanin magic
recipe \cite{Sk95}, to each of the Lax matrices. Actually, it has
also been shown that integration of the Hamilton--Jacobi equation
involves Abelian dif\/ferentials of order independent of the number of
sites. Moreover, for $\alg=sl(r)$ the genus of the underlying
Riemann surface can be computed to be $g=\frac{r(r-1)}{2}$.
\end{itemize}
For further use, we remark that the ring of spectral invariants of a
Lax matrix is, as it is well known, invariant w.r.t.\ rational
changes in the spectral parameter and/or multiplication by a~rational function of the spectral parameter. Thus, instead of the
Hamiltonians \eqref{ham-fm} associated with the ``parameter free'' Lax
matrices \eqref{lax-fm}, we can consider, as an equivalent set of
Hamiltonians for the generalised $gl(r)$ bending f\/lows, the spectral
invariants associated with the  matrices
\begin{gather}
L_1=\frac{X_2}{z-z_2}+\frac{X_1}{z-z_1},\qquad
L_2=\frac{X_3}{z-z_2}+\frac{X_1+X_2}{z-z_1}, \qquad \ldots,\nonumber\\
L_{N-1}=\frac{X_N}{z-z_2}+\frac{\sum\limits_{i=1}^{N-1} X_i}{z-z_1},\label{ratlaxmat}
\end{gather}
rationally dependent on $z$. We f\/inally want to add the following
observation. As it is pictorially suggested in Fig.~\ref{fig2.4},
there is a natural interpretation of our procedure in terms of
triangulations of a polygon and graphs.

Actually, we have already seen (see Fig.~\ref{figurgr}) that we can
interpret our procedure of glueing points in terms of building up tree-graphs
out of the set of singular points of the generic Lax matrix of the Gaudin
model. In particular, the glueing pattern depicted on the right of Fig.~\ref{fig2.4} corresponds to the limit in which we produce the
``Bending'' integrals, and, in the topological point of view, exactly gives
(provided we identify the point~$z_1$ with the root), the dual graph to the
triangulation of the right hand side of the picture.

It is outside of the size of the paper to fully discuss these
issues, which however the authors think being of a certain interest.
It should also be noticed\footnote{We thank G.~Felder for an
observation concerning this point.} that tree diagrams corresponding
to the glueing patterns described in this paper appear in the theory
of moduli of pointed rational curves.

\section{The quantum case}

\subsection{FFR Gaudin algebras and their limits}
The existence of a large quantum commutative subalgebra
$\cA(z_1,\dots,z_N)\subset U(\alg)^{\otimes N}$ containing the
Gaudin Hamiltonians $H_i$ was shown in \cite{FFR94}.  As we recalled
in the introduction, whenever $\alg$ has rank greater than one,
$\cA(z_1,\dots,z_N)$ has new generators, besides the quadratic ones,
known as {\em higher Gaudin Hamiltonians}.

Roughly speaking, the def\/inition  of $\cA(z_1,\dots,z_N)$ can be
obtained identifying a commutative subalgebra of the enveloping
algebra $U(\gg\otimes t^{-1}\mathbb{C}[t^{-1}])$ as follows. To any
collection $z_1,\dots,z_n$ of pairwise distinct complex numbers, one
can naturally assign the evaluation map $U(\alg\otimes
t^{-1}\mathbb{C}[t^{-1}])\to U(\gg)^{\otimes N}$. The image of the
commutative subalgebra in $U(t^{-1}\gg[t^{-1}])$
under the composition of the above homomorphisms, denoted by
$\cA(z_1,\dots,z_n)$,  is called (\textit{quantum}) \emph{Gaudin algebra}.

It is useful to recall the following well known facts from the
theory of Lie algebras. Let $\alg$ be a Lie algebra. The universal
enveloping algebra $U(\gg)$ bears a natural f\/iltration by the degree
with respect to the generators. The associated graded algebra is the
symmetric algebra $S(\gg)$ by the Poincar\'e--Birkhof\/f--Witt
theorem. To every element $\xi$ in $U(\gg)$ there corresponds
uniquely its image $\bar{\xi}$ in $S(\gg)$, since $S(\gg)$
corresponds to the graded algebra associated with $U(\gg)$. This
element is customarily called, in the theory of Gaudin and
Heisenberg spin systems/chains, the ``classical limit'' of $\xi$.
Also, whenever $\alg$ is reductive, $S(\gg)\simeq S(\gg^*)$, and
$\dsl{ U(\underbrace{\alg\oplus\alg\cdots\oplus\alg}_{N\,{\rm
times}})=\big(U(\alg)^{\otimes^N}\big)}.$

The problem of f\/inding explicit representatives for generators of
$\cA(z_1,\dots,z_N)$ in the case $\alg=gl(r)$ was concretely solved
by Talalaev some ten years later after the paper by Feigin, Frenkel
and Reshetikhin \cite{FFR94}, by means of the following
construction.
\begin{theorem}[Talalaev, 2004]\label{tathe04}
Let  $L(z)$ be the Lax matrix of the $gl(r)$-Gaudin model, that is,
let $L(z)$  satisfy the $r$-matrix commutation relations
\begin{equation}\label{rmatcommrel}
[ L(z)\otimes 1, 1 \otimes L(u)] = \left[\frac{\Pi}{z-u}, L(z)\otimes 1 +
1 \otimes L(u)\right]. \end{equation}
Consider the differential operator
in the variable $z$ $\text{\rm ``DET''}\big(\partial_z - L(z)\big)=
\sum\limits_{i=0,\dots,r} QH_{i}(z)
\partial_z^i$;
Then:
\begin{equation*}
\forall \, i,j\in 0,\ldots,r, \quad \text{and} \quad u,v\in \CC, \qquad
[QH_{i}(z)|_{z=u} , QH_{j}(z)|_{z=v} ]=0.
\end{equation*}
\end{theorem}
The $QH_i(z)$'s generators of (a full set of) quantum mutually
commuting quantities. We shall discuss the meaning of the symbol
$\text{``DET''}$ appearing in the formulation of Talalaev's theorem in
Section \ref{ManinM}, where we shall also comment on its proof.

For the moment, let us remark that the Gaudin algebra(s) def\/ined by
Feigin, Frenkel and Reshetikhin -- and concretely identif\/ied by
Talalaev's formula for $N$-spin $\mathfrak{gl}_r$ Gaudin systems~--
explicitly depends, in general (that is when $N\ge 3$)
on the points $z_1, \ldots,z_N$.
In another parlance, its has moduli. However it holds
(see \cite{CFR1} for the proof) the following
\begin{proposition}
$\cA(z_1,\ldots,z_N)$ is invariant under permutations and under
simultaneous re\-sca\-lings $ z_i\to \alpha z_i+\beta$, $\alpha,\beta\in
\cc $; thus the ``two site'' algebra $\cA(z_1,z_2)$ is independent of
$z_1$, $z_2$.
\end{proposition}

Let us now consider limits of the quantum Gaudin algebras when some
of the points $z_1,\dots,z_N$ glue together according to some
pattern. Here, obviously enough, we will follow the same glueing
procedure as in Section \ref{sect-class}.

We recall that the basic procedure is that we  keep some (say the
f\/irst $k$) points $z_1,\dots,z_k$ ``f\/ixed'', and let the remaining
$N-k$ points glue to a new point $w$, via
\begin{equation*}
z_{k+i}=w+s u_i,\qquad  i=1,\dots,N-k,\quad z_i\neq z_j,\quad u_i\neq
u_j, \quad s\to 0 .
\end{equation*}
To describe this limit in the quantum case, it is useful to
introduce the maps
\[ D_{k,N}:=\text{id}^{\otimes
k}\otimes\diag_{N-k}: \ U(\gg)^{\otimes (k+1)}\hookrightarrow
U(\gg)^{\otimes N},
\]
 def\/ined by
 \[
D_{k,N}(X_1\otimes\cdots\otimes X_k\otimes X_{k+1})=X_1\otimes\cdots
X_{k}\otimes\underbrace{X_{k+1}\otimes\cdots\otimes X_{k+1}}_{N-k
\text{ times}}
\] and
\[
I_{k,N}:=1^{\otimes k}\otimes\text{id}^{\otimes(N-k)}: \
U(\gg)^{\otimes (N-k)}\hookrightarrow U(\gg)^{\otimes N},
\]
 def\/ined by
 \[
  I_{k,N}(X_1\otimes\cdots
X_{N-k})=\underbrace{\mathbf{1}\otimes\cdots\otimes\mathbf{1}}_{k\text{
times}}\otimes X_1\otimes\cdots\otimes X_{N-k}.
\]
Furthermore, let us def\/ine the algebra
\[
\A_{(z_1,\dots,z_k,w),(u_1,\dots,u_{N-k})}
:=D_{k,N}(\A(z_1,\dots,z_k,w))\cdot I_{k,N}(\A(u_1,\dots,u_{N-k}))
\]
In \cite{CFR1} it is proven the following

\begin{theorem}
The algebra $\A_{(z_1,\dots,z_k,w),(u_1,\dots,u_{N-k})}$ is
commutative;
\[\lim\limits_{s\to 0}\A(z_1,\dots,z_k,w+su_1,\dots,z+s\,u_{N-k})=
\A_{(z_1,\dots,z_k,z),(u_1,\dots,u_{N-k})}.\]
\end{theorem}

\noindent
{\bf Remarks.} 1) In some sense we arrive at a kind of factorisation
of the limit algebra by ``adding'' one point. Indeed, the commutative
algebra $\A_{(z_1,\dots,z_k,w),(u_1,\dots,u_{N-k})}$ involves the
FFR algebra associated with the points $(z_1,\ldots,z_k,w)$, and a
FFR algebra related with the points $u_1,\ldots,u_{N-k}$.
Iterating this limiting procedure described above we can obtain some
new commutative subalgebras in $U(\gg)^{\otimes N}$, just like in
the classical case we found Poisson commutative subalgebras of~$S(\alg)^{\otimes^N}$.

2) In the case $\alg=gl(r)$ the classical limits (in the sense
specif\/ied above) of, respectively, $I_{k,N}(\A(u_1,\dots,u_{N-k}))$ and
$D_{k,N}(\A(z_1,\dots,z_k,w))$ coincide with the ring of the spectral
invariants of the matrices
\[
L_1= \sum\limits_{i=k+1}^N\frac{X_i}{z-u_i}  , \qquad {\rm and}\qquad
 L_2(w)= \sum\limits_{i=1}^k \frac{X_1}{z-z_i}+\frac{\sum\limits_{m=k+1}^N
X_m}{z-w}
\] already discussed in Section \ref{sect-class}.

In particular, as in the bending f\/low case, we can iterate the
procedure to pass from $\cA(z_1,z_2$, $\ldots,z_N)$ to $
\cA_{(z_1,z_2,\ldots,z_{N-2},w)\,(z_{N-1},z_N)}$, and f\/inally obtain
the subalgebra
\begin{equation}\label{Alim}
\A_{\lim}\equiv\A_{(z_1,z_2),\dots,(z_1,z_2)}\subset U(\gg)^{\otimes
n}.
\end{equation}
In the process, we have to use translation invariance of
$\cA(z_1,\ldots,z_M)$, as well as the property that
the two site Gaudin algebra $\cA(u,v)$ is
independent of $(u,v)$.

The limit algebra \eqref{Alim} is generated by
\begin{equation*}
D_{1,N}(\A(z_1,z_2)),\quad 1\otimes
D_{1,N-1}(\A(z_1,z_2)),\quad \dots,\quad 1^{\otimes(N-2)}\otimes
\A(z_1,z_2).
\end{equation*}
More explicitly:

\begin{proposition}
The subalgebra $\A_{\lim}$ is generated by elements
$H_{l,k}^{(\alpha)}\in U(gl(r))^{\otimes n}$ such that their
classical ``limits'' $\overline H_{l,k}^{(\alpha)}$, are given by
\[
\overline
H_{l,k}^{(\alpha)}(X_1,\dots,X_n):=\mathop{\rm Res}\limits_{z=0}\frac{1}{z^{\alpha+1}}
\mathop{\rm Tr}\nolimits\left(X_k+z \sum\limits_{i=k+1}^nX_i\right)^l,
\]
\end{proposition}
In other words, the classical limits of the (still unspecif\/ied)
``quantum Bending Hamiltonians'' do coincide with the spectral
invariants of the Lax matrices $ L_k(z)=z X_k+\sum\limits_{i=k+1}^N X_i$,
$i=2,\ldots,N$, associated with the bending f\/lows. In turn, these can
be traded for the spectral invariants of their ``rational'' analogs
\begin{gather}
L_1(z)=\frac{X_{N-1}}{z-z_2}+\frac{X_N}{z-z_1},\qquad
L_2(z)=\frac{X_{N-2}}{z-z_2}+\frac{X_N+X_{N-1}}{z-z_1},\qquad \ldots,\nonumber\\
L_{N-1}=\frac{X_1}{z-z_2}+\frac{\sum\limits_{i=2}^{N} X_i}{z-z_1}.\label{cqlax-bend}
\end{gather}

In the last part of the paper we shall show how, using Talalaev's
results brief\/ly reminded before as well as some of the results of~\cite{CF07} we can f\/ind a suitable quantization of the traces of the
powers of the Lax matrices~\eqref{cqlax-bend} (as well as of, more
generally, of the quantum Lax matri\-ces~$L_1(z)$,~$L_2(z)$ associated
with an elementary  ``gluing'' procedure.

To this end, we need to introduce the notion of Manin matrix.

\subsection{Manin matrices and the quantization of traces of Lax matrices}\label{ManinM}

We consider the (dif\/ferential operator valued) matrix considered in
Talalaev's Theorem~\ref{tathe04}, that is, $ \mathcal{M}(z)=\partial_z-L(z)$. The
standard linear $r$-matrix commutation relations~\eqref{rmatcommrel}
imply that the matrix elements of $\mathcal{M}$ satisfy special
commutation relations. These are the def\/ining commutation relations
of a class of matrices (with non-commutative entries) called in~\cite{CF07} {\em Manin matrices.} The terminology originates from a
seminal paper of Yu.I.~Manin on quantum groups~\cite{Ma88}.

Manin matrices are, in a suitable sense, matrices associated with
linear maps between commutative rings. Their (more operative)
def\/inition can be given as follows:
\begin{definition}
 Let $M_{ij}$ be a matrix with
elements in a non commutative (unital) ring $\mathcal{R}$; we call
it a (column) Manin matrix if:
\begin{itemize}\itemsep=0pt
\item  elements in the same column commute among themselves;
\item  commutators of the cross terms in any $2\times 2$ submatrix
are equal: \[ [M_{ij}, M_{kl}]=[M_{kj}, M_{il}],\qquad \text{e.g.} \quad
[M_{11}, M_{22}]=[M_{21}, M_{12}] \quad \text{and so on and so forth}.\]
\end{itemize}
\end{definition}

Manin matrices admit a natural def\/inition of determinant. Indeed, if
$M$ is {\em any} matrix, one can in principle def\/ine ``a''
determinant of $M$ by column expansion:
\begin{equation*}
\det M=\det\nolimits^{\rm col} M=\sum_{\sigma\in
S_n} (-1)^\sigma \prod^{\curvearrowright}_{i=1,\dots,n}
M_{\sigma(i),i}, \end{equation*} where $S_n$ is the group of
permutations of $n$ letters, and the symbol $\curvearrowright$ means
that in the product $\prod\limits_{i=1,\dots,n} M_{\sigma(i),i}$ one writes
at f\/irst the elements from the f\/irst column, then from the second
column and so on and so forth.

In the case of \MMs\ the characteristic property of total
antisymmetry of the determinant is preserved:
\begin{proposition} The column determinant of a \MM~ does not depend on the order of the
columns in the column expansion, i.e.,
\begin{gather*}
 \forall \, p\in S_n \qquad
\det\nolimits^{\rm col} M=\sum_{\sigma\in S_n} (-1)^\sigma
\prod^{\curvearrowright}_{i=1,\dots,n} M_{\sigma(p(i)),p(i)}.
\end{gather*}
\end{proposition}
In particular, a  good notion of ``DET'' in the formulation of
Talalaev's quantum integrals of the motion is given by the column
determinant of~$\mathcal{M}(z)$.

Among the properties of ``ordinary'' linear algebra that \MMs\
satisfy, we mention the following:
\begin{itemize}\itemsep=0pt
\item
The inverse of a \MM~ $M$, whenever it exists, is still a  Manin matrix.
\item The (left) Cramer formula holds:
\[
M^{\rm adj} M= \det\nolimits^{\rm col}(M) \mathbf{1},
\]
where $M^{\rm adj}$ is computed -- via column determinants -- in the
same way as in the ordinary case.
\item Schur's formula for the determinant of block matrices holds:
\begin{gather*}
 \det\nolimits^{\rm col}\left(\begin{array}{cc}
A & B\\
C & D \\
\end{array}\right)=\det\nolimits^{\rm col}(A)\det\nolimits^{\rm col}(D-CA^{-1}B)=\det\nolimits^{\rm col}(D)\det\nolimits^{\rm col}
\big(A-BD^{-1}C\big) .
\end{gather*}
\item
The Cayley--Hamilton theorem, $\det\nolimits^{\rm col}(t-M)|_{t=M}=0$, holds true.
\end{itemize}

Further properties of \MMs, as well as full proofs of these facts
can be found in~\cite{CFR08}. Particularly important, for our
purposes, is the the fact that a suitable redef\/inition of the traces
of the powers of the quantum Lax matrix allows us to trade
coef\/f\/icients of the characteristic polynomials for these ``normalised
traces''. The key property concerns the quantum counterparts of the
classical Newton identities between these two classes of
${\rm Ad}$-invariant quantities. The coef\/f\/icients $\sigma_i$ of the
characteristic polynomial of a $n\times n$ matrix $M$ are, as it is
well known, the elementary symmetric polynomials in the eigenvalues
$\la_i$ of $M$; however in the ring of symmetric polynomials in the
$n$ variables $\la_i$, one can consider other sets of generators and
in particular the so called power sums
\begin{equation*}
    \tau_i=\sum_{k=1}^{n} \lambda_{k}^{i}=\text{Tr}\,\big(M^i\big), \qquad i=1,
\ldots,n  .
\end{equation*}
In the case of matrices with commuting entries, the family
$\{\sigma_i\}$, ${i=1,\ldots,n}$ and $\{\tau_i\}$, $i={1,\ldots,n}$ are
related by the identities (the {\em Newton} identities):
\begin{equation}\label{newid}
    (-1)^{k+1} k    \sigma_k =\sum_{i= 0,\dots,k-1}  (-1)^{i} \sigma_{i}
    \tau_{k-i}.
\end{equation}
The same identities hold true, in the case of \MMs, for the
corresponding quantities.
\begin{theorem}
  Newton identities of the form \eqref{newid} between $\tau_i={\rm Tr}\,M^k$ and the
  coefficients $\sigma_i$ of the expansion of $\det\nolimits^{\rm col}(t+M)$ in powers of $t$ hold.
\end{theorem}

\begin{proof}
First one can observe that the following property, whose proof is
straightforward,
\[
{\rm Tr}\, (t+M)^{\rm adj}=\partial_t \det\nolimits^{\rm col}(t+M)
\]
hold. Then we have that
\begin{gather*}
 \frac{1}{t} \sum_{k=0,\dots,\infty}
{\rm Tr}\,\big((-M/t)^k\big) = {\rm Tr}\,
\frac{1}{t+M}={\rm Tr}\, \big(\big(\det\nolimits^{\rm col}(t+M)\big)^{-1}(t+M)^{\rm adj} \big)  \\
\qquad{} = \big(\det\nolimits^{\rm col}(t+M)\big)^{-1} {\rm Tr}\, (t+M)^{\rm adj}=
\big(\det\nolimits^{\rm col}(t+M)\big)^{-1}
\partial_t \det\nolimits^{\rm col}(t+M).
\end{gather*} The result is obtained substituting $-M$ for $M$ in these
formulas.
\end{proof}

The quantization of the traces of powers of the Lax matrices can be
def\/ined according to the Newton identities and the following
considerations.
First, a straightforward computation that makes use of the the
$r$-matrix commutation relations~\eqref{rmatcommrel}, proves the
following
\begin{proposition}
Let\,\footnote{In these and subsequent formul\ae\ we use the symbol
$\hat{}$ to stress the fact that in the quantum case $\hat{X_i}$ are
to be regarded as elements in the universal enveloping algebra of
$gl(r)$.}
\begin{equation}\label{quantgaulax}
\hat{L}=\sum_{i=1}^N\dsl{\frac{\hat{X}^i}{z-z_i}}
\end{equation}
be the Lax matrix for the quantum Gaudin model.
 For any number of sites $N$, and different $z_i$ the
matrix $\partial_z - \hat{L}(z)$ is a~\MM.
\end{proposition}

It might be  mentioned that this is not a specif\/ic property  for the
Gaudin system herewith considered. Indeed, Manin matrices enter
other topics in (quantum) integrability; for instance,
$e^{-{\partial_z}} T_{gl(r)\text{-Yangian}}(z)$ is Manin, where
$T_{gl(r)\text{-Yangian}}(z)$ is the Lax (or ``transfer'') matrix for the
Yangian algebra $Y(gl(n))$, satisfying quadratic $r$-matrix
commutation relations (see~\cite{CF07} for more details).

The ``quantization'' of the traces of a Lax matrix of the form~\eqref{quantgaulax}  can be achieved as follows. From Talalaev's
theorem we know that the coef\/f\/icients $QH(z)$ of characteristic
polynomial $\det\nolimits^{\rm col}(\partial_z-\hat{ L}(z))$ do commute (at the quantum
level), and from the theorem about Newton identities we know that we
can trade these coef\/f\/icients for traces of the powers of the
corresponding Manin matrix.

Basically, we simply have to remark that we should not consider the
quantities $\dsl{\text{Tr}(\hat{L}(z)^k})$, which, as it has been
shown in \cite{CRT} do not, generally speaking, commute, but rather
the traces of powers of the Manin matrix
\begin{equation}\label{quantrace}
\text{Tr}\big((\partial_z-\hat{L}(z))^k\big)=
\sum_{j=0}^k (Q\,{\rm Tr})^{k}_j(z)\partial_z^{k-j},\qquad k=1,\ldots, r.
\end{equation}
As it is easily seen, there is a recursion relation of the form
$  Q\,{\rm Tr}^{k}_{j+1}(z)\simeq Q\,{\rm Tr}^{k+1}_j(z) $, and  hence, to
obtain the expected number of independent quantities, we can
consider simply the coef\/f\/icients $Q\,{\rm Tr}^{k}_{k}$, that is the
coef\/f\/icients of zeroth order of each dif\/ferential ``polynomial''~\eqref{quantrace}.

In turn, one sees that these quantities are given by the traces of
matrices  ${\hat L}_k(z)^{[n]}$, that can be called "quantum powers" of
$L(z)$;
these powers are
def\/ined by the Fa\`a di Bruno formula (see, e.g.,~\cite{Di})
\[
{\hat L}_k^{[0]}(z) = {\rm Id},\qquad {\hat L}_k^{[i]}(z) = {\hat
L}_k^{[i-1]}(z) {\hat L}_k(z) - \frac{\partial}{\partial z} \big({\hat
L}_k^{[i-1]}(z)\big).\]

These arguments work for any Lax matrix of ``Gaudin type'', that is,
of the form~\eqref{quantgaulax}; in particular, they hold for the
(rational) Lax matrices $L_i(z)$ \eqref{ratlaxmat} associated with
the Bending f\/lows, since the dependence on the spectral parameter is
chosen in such a way to fulf\/il the canonical $r$-matrix commutation
relations \eqref{rmatcommrel}. Thus computing the traces of the
quantum powers of each of the (quantum counterparts of the) matrices
$L_i(z)$ of~\eqref{ratlaxmat} one gets a quantization of the
Hamiltonians of the Bending f\/lows.

\section{Concluding remarks and open problems}
In this paper we have discussed, both in the classical and in the
quantum settings, some features in the theory of (homogeneous)
Gaudin models, and concentrated on the  Gaudin algebras of mutually
commuting quantities. We have addressed, in particular, the problem
of a suitable def\/inition of a complete set of Hamiltonians in the
case when the poles of the Lax matrix $L_G(z)$ of the model glue
together. In the classical case, we have shown that a complete set
of mutually commuting Liouville integrals can be obtained by a
suitable analysis of the limiting procedure on $L_G(z)$. We have then
addressed some features of the (multi)-Hamiltonian properties of the
limits of these (classical) Gaudin systems, in the framework of the
Hamiltonian geometry of Loop algebras. Also, we touched the problem
of framing the geometry of these systems within the combinatorial
approach to the geometry of moduli spaces of pointed rational
curves. In these respects, the results herewith collected are
somewhat sketchy, and will be developed elsewhere.

In the quantum framework, we have followed the line of \cite{FFR94}
as well as previous results of ours~\cite{CFR1,FM06}. We reviewed how our
gluing procedure for the poles of the Lax matrix of the Gaudin model
can be applied in the quantum case to yield new quantum Gaudin
algebras. Also, by means of results of~\cite{CF07} about the theory
of the so-called Manin matrices, we have shown quite explicitly the
procedure to def\/ine, in the quantum Gaudin algebras, a set of
generators that reduce in the classical case to the usual traces of
powers of the Lax matrix.

Besides the full solution of the problems in the Hamiltonian and
combinatorial aspects of these systems that were just addressed in
this paper, we can envisage at least two domains where the study of
the limits discussed in the present paper might be of some interest.
The f\/irst is the realm of Heisenberg-like models, and the theory of
Yangian algebras. The second is the study of non-homogeneous Gaudin
models, that, besides of their own interest in the mathematical
theory of integrable (quantum) systems, proved to be important in
the theory of strongly correlated electron systems (see, e.g.,
\cite{Si02}). Due to the fact that non homogeneous systems lack
global $G$-invariance, however, the results of this paper cannot be
applied to that case in a straightforward manner. Work in this
direction is in progress.

\subsection*{Acknowledgments}

 G.F. acknowledges support from the ESF
programme MISGAM, and the Marie Curie RTN ENIGMA.
The work of A.C. and L.R. has been
partially supported
by the Russian Federal Nuclear Energy Agency. A.C. has been
partially supported by the Russian President Grant MK-5056.2007.1,
by the grant of Support for the Scientif\/ic Schools 3035.2008.2, by
the RFBR grant 08-02-00287a, and by the ANR grant GIMP (Geometry and
Integrability in Mathematics and Physics). The work of L.R. was
partially supported by the RFBR grant 05 01 00988-a and the RFBR
grants 05-01-02805-CNRSL-a and 07-01-92214-CNRSL-a. Also, L.R.
gratefully acknowledges the support from the Deligne 2004 Balzan
prize in mathematics. We thank the referees for useful remarks.

\pdfbookmark[1]{References}{ref}
\LastPageEnding


\begin{thebibliography}{99}

\footnotesize\itemsep=0pt



\bibitem{BBT03} Babelon O., Bernard D., Talon M.,  Introduction to classical integrable systems,
{\it Cambridge Monographs on Mathematical Physics}, Cambridge University Press, Cambridge, 2003.



\bibitem{Ch03} Chernyakov Yu.B.,
 Integrable systems, obtained by point fusion from rational and
elliptic Gaudin systems, {\it Theoret. and Math. Phys.} {\bf 141} (2004), 1361--1380,
\href{http://arxiv.org/abs/hep-th/0311027}{hep-th/0311027}.

\bibitem{CF07}
 Chervov A., Falqui G.,  Manin  matrices and Talalaev's
formula, {\it J. Phys. A: Math. Theor.} {\bf 41} (2008), 194006, 28~pages,
\href{http://arxiv.org/abs/0711.2236}{arXiv:0711.2236}.

\bibitem{CFR1} Chervov A.,  Falqui G., Rybnikov L.,  Limits of Gaudin algebras, quantization of bending f\/lows, Jucys--Murphy elements and
Gelfand--Tsetlin bases, \href{http://arxiv.org/abs/0710.4971}{arXiv:0710.4971}.

\bibitem{CFR08} Chervov A., Falqui G., Rubtsov V.,
Algebraic properties
    of Manin matrices I, {\it Adv. Appl. Math.}, to appear, \href{http://arxiv.org/abs/0710.4971}{arXiv:0901.0235}.

\bibitem{CRT} Chervov A., Rybnikov L., Talalaev D.,
Rational Lax operators and their quantization,
\href{http://arxiv.org/abs/hep-th/0404106}{hep-th/0404106}.

\bibitem{CT04} Chervov A., Talalaev D.,
Universal G-oper and Gaudin eigenproblem,
\href{http://arxiv.org/abs/hep-th/0409007}{hep-th/0409007}.

\bibitem{CT06}  Chervov A., Talalaev D.,
Quantum spectral curves, quantum integrable systems and the
geometric Langlands correspondence,
\href{http://arxiv.org/abs/hep-th/0604128}{hep-th/0604128}.

\bibitem{Di}
 Dickey L., Soliton equations and Hamiltonian systems, 2nd ed.,
{\it Advanced Series in Mathematical Physics}, Vol.~26,
World Scientif\/ic Publishing Co., Inc., River Edge, NJ, 2003.

\bibitem{ER96} Enriquez B., Rubtsov V.,  Hitchin systems,
higher Gaudin operators and $R$-matrices,
{\it Math. Res. Lett.}  {\bf 3}   (1996), 343--357,
\href{http://arxiv.org/abs/alg-geom/9503010}{alg-geom/9503010}.

\bibitem{FM01}
 Falqui G., Musso F.,  Gaudin models and bending f\/lows: a
geometrical point of view, {\it J. Phys. A: Math. Gen.} {\bf 36}  (2003), 11655--11676,
\href{http://arxiv.org/abs/nlin.SI/0306005}{nlin.SI/0306005}.

\bibitem{FM03-2}
 Falqui G., Musso F.,  Bi-Hamiltonian geometry and separation of
variables for Gaudin models: a case study, in SPT 2002: Symmetry and Perturbation Theory (Cala Gonone), Editors S.~Abenda, G.~Gaeta and S.~Walcher, World Sci. Publ., River Edge, NJ, 2002, 42--50, \href{http://arxiv.org/abs/nlin.SI/0306008}{nlin.SI/0306008}.

\bibitem{FM04}
 Falqui G., Musso F., On separation of variables for homogeneous
${\rm sl}(r)$ Gaudin systems,
 {\it Math. Phys. Anal. Geom.}   {\bf 9}  (2006), 233--262,
\href{http://arxiv.org/abs/nlin.SI/0402026}{nlin.SI/0402026}.

\bibitem{FM06}
 Falqui G., Musso F.,  Quantisation of bending f\/lows,
 {\it Czechoslovak J. Phys.}   {\bf 56}  (2006), 1143--1148,
\href{http://arxiv.org/abs/nlin.SI/0610003}{nlin.SI/0610003}. 

\bibitem{FaPe03}  Falqui G., Pedroni M.,  Separation
of variables for bi-Hamiltonian systems, {\em Math. Phys. Anal.
Geom.} {\bf 6}  (2003), 139--179, \href{http://arxiv.org/abs/nlin.SI/0204029}{nlin.SI/0204029}.

\bibitem{FlM01}
 Flaschka H., Millson J.,  Bending f\/lows for sums of rank one matrices,
{\it Canad. J. Math.}  {\bf 57},  (2005), 114--158,
\href{http://arxiv.org/abs/math.SG/0108191}{math.SG/0108191}. 

 \bibitem{FF} Feigin B., Frenkel E., Af\/f\/ine Kac--Moody
algebras at the critical level and Gelfand--Dikii algebras,
in Inf\/inite Analysis, Part~A,~B (Kyoto, 1991), {\it Adv. Ser. Math. Phys.}, Vol.~16, World Sci. Publ., River Edge, NJ, 1992, 197--215.

\bibitem{FFR94} Feigin B., Frenkel E., Reshetikhin N.,  Gaudin
    model, Bethe ansatz and critical level, {\it Comm. Math. Phys.} {\bf 166}
  (1994), 27--62, \href{http://arxiv.org/abs/hep-th/9402022}{hep-th/9402022}.

\bibitem{FFTL} Feigin B., Frenkel E., Toledano-Laredo V.,  Gaudin model with irregular singularities,
\href{http://arxiv.org/abs/math.QA/0612798}{math.QA/0612798}. 

\bibitem{Fr95}
 Frenkel E.,
Af\/f\/ine algebras, Langlands duality and Bethe ansatz, in Proceedings XIth
International Congress of Mathematical Physics (Paris, 1994), Int. Press, Cambridge, MA, 1995,
606--642, \href{http://arxiv.org/abs/q-alg/9506003}{q-alg/9506003}.

\bibitem{G83} Gaudin M., La fonction d'onde de Bethe,
Collection du Commissariat \`a l'\'Energie
Atomique: S\'erie Scientif\/ique,
Masson, Paris, 1983. 

\bibitem{GZ00}
 Gel'fand I.M., Zakharevich I.,  Webs, Lenard schemes, and the
local geometry of bi-Hamiltonian Toda and Lax structures,
{\it Selecta Math. (N.S.)} {\bf 6} (2000), 131--183, \href{http://arxiv.org/abs/math.DG/9903080}{math.DG/9903080}.

\bibitem{Ju89} Jur\v{c}o B.,  Classical Yang--Baxter equations and
quantum integrable systems, {\it J. Math. Phys.} {\bf 30} (1989), 1289--1293.

\bibitem{KM96}
 Kapovich M., Millson J.,   The symplectic geometry of polygons
in Euclidean space, {\it J. Differential Geom.} {\bf 44} (1996), 479--513.



\bibitem{Ku92}
 Kuznetsov V.B., Quadrics on Riemannian spaces of constant
curvature. Separation of variables and a connection with the Gaudin
magnet, {\it Theoret. and Math. Phys.}  {\bf 91} (1992), 385--404.

\bibitem{Ma88}
  Manin Yu.I.,  Quantum groups and noncommutative geometry,
Universit\'e de Montr\'eal, Centre de Recherches Math\'ematiques, Montreal, QC, 1988.

\bibitem{MPR05}
 Musso F., Petrera M., Ragnisco O.,  Algebraic extensions of
Gaudin models,  {\it J. Nonlinear Math. Phys.} {\bf  12}  (2005), suppl.~1,
482--498, \href{http://arxiv.org/abs/nlin.SI/0410016}{nlin.SI/0410016}.

\bibitem{PV} Pedroni M., Vanhaecke P.,
A Lie algebraic generalization of the Mumford system,
its symmetries and its multi-Hamiltonian structure,  {\it Regul. Chaotic
Dyn.} {\bf 3}  (1998), 132--160.

\bibitem{RSTS} Reyman A.G., Semenov-Tian-Shansky M.A., Group-theoretical methods in the theory of f\/inite-dimensional integrable systems, in Dynamical Systems VII, Editors V.I.~Arnold and S.P.~Novikov, {\it Encyclopaedia of Mathematical Sciences}, Vol.~16, Berlin, Springer, 1994, 116--225.

\bibitem{Ryb1} Rybnikov L.G., The shift of invariants method and the Gaudin model, {\it Func. Anal. Appl.} {\bf 40} (2006), 188--199,
\href{http://arxiv.org/abs/math.RT/0606380}{math.RT/0606380}. 


\bibitem{Ryb}  Rybnikov L.G., Uniqueness
of higher Gaudin Hamiltonians,
\href{http://arxiv.org/abs/math.QA/0608588}{math.QA/0608588}. 

\bibitem{Si02} Sierra G.,   Integrability and conformal symmetry in BCS model,
in Statistical Field Theories (Como, 2001), {\it  NATO Sci. Ser. II Math. Phys. Chem.}, Vol.~73, Kluwer Acad. Publ., Dordrecht, 2002, 317--328, \mbox{\href{http://arxiv.org/abs/hep-th/0111114}{hep-th/0111114}}.

\bibitem{Sk95} Sklyanin E.K.,  Separation of variables~--
new trends, in
Quantum Field Theory, Integrable Models and Beyond (Kyoto, 1994),
{\it Progr. Theoret. Phys. Suppl.}  (1995), no.~118, 35--60, \href{http://arxiv.org/abs/solv-int/9504001}{solv-int/9504001}.

\bibitem{Ta04} Talalaev D.,  Quantization of the Gaudin system,
{\it Funct. Anal. Appl.} {\bf 40} (2006), 73--77,
\href{http://arxiv.org/abs/hep-th/0404153}{hep-th/0404153}.

\end{thebibliography}
\end{document}